\documentclass[a4paper, 10pt]{article} % Font size (can be 10pt, 11pt or 12pt) and paper size (remove a4paper for US letter paper)

\usepackage{graphicx}
\usepackage{amsmath,amssymb,amsfonts,amsthm,latexsym}
\usepackage{mathrsfs}
\usepackage{textcomp}
\usepackage{tikz-cd}
\usepackage[all]{xy}

\usepackage{mathpazo} % Use the Palatino font
\newtheorem{thm}{Theorem}
\newtheorem{cor}[thm]{Corollary}
\newtheorem{lem}[thm]{Lemma}
\newtheorem{prop}[thm]{Proposition}
\theoremstyle{definition}
\newtheorem{defn}[thm]{Definition}
\theoremstyle{remark}
\newtheorem{rem}[thm]{Remark}
\theoremstyle{remark}
\newtheorem{ex}[thm]{Example}
\theoremstyle{remark}

\theoremstyle{remark}

\newcommand{\N}{\mbox{$\mathbb N$}}
\newcommand{\R}{\mbox{$\mathbb R$}}

\newcommand\set[1]{\mbox{$\{#1\}$}}
\newcommand\cat[1]{\mbox{$\mathbf{#1}$}}
\newcommand\ds[1]{\mbox{$\displaystyle{#1}$}}
\newcommand\map[3]{\mbox{$#1:#2\rightarrow #3$}}
\newcommand\inv[1]{\mbox{$#1^{-1}$}}
\newcommand\trip[1]{\mbox{$\mathbb{#1}$}}
\newcommand\lift[1]{\mbox{$\mathfrak{#1}$}}

\newcommand{\id}{\mbox{$\mathfrak{I}$}}
\newcommand{\eps}{\mbox{$\ds{\epsilon}$}}

\newcommand\comp[2]{\mbox{$#1\cdot #2$}}
\newcommand\ocomp[2]{\mbox{$#1\circ #2$}}
\newcommand\join[2]{\mbox{$#1\vee #2$}}
\newcommand\meet[2]{\mbox{$#1\wedge #2$}}
\newcommand\mon[3]{\mbox{$(#1,#2,#3)$}}
\newcommand\powcat[2]{\mbox{$#1^{\mathbb{#2}}$}}
\newcommand\alg[1]{\mbox{$\mathbf{Alg}\left(\mathbb{#1}\right)$}}

\newcommand\Adj[4]{\mbox{$$\xymatrix{#1 \ar@<1ex>[rrr]^{#2}_*-<0.55ex>{_\perp} & & & \ar@<1ex>[lll]^{#3} #4}$$}}
\newcommand\compAdj[7]{\mbox{$$\xymatrix{#1 \ar@<1ex>[rrr]^{#2}_*-<0.55ex>{_\perp} & & & \ar@<1ex>[lll]^{#3} #4 \ar@<1ex>[rrr]^{#5}_*-<0.55ex>{_\perp} & & & \ar@<1ex>[lll]^{#6} #7}$$}}

\title{\textbf{A pair of monads in Topology}\\ % Title
} % Subtitle

\author{\textsc{A. Razafindrakoto\thanks{Email: arazafindrakoto@uwc.ac.za - Department of Mathematics and Applied Mathematics - University of the Western Cape}}} % Institution

\date{} % Date

\begin{document}

\maketitle % Print the title section

%----------------------------------------------------------------------------------------
%	ABSTRACT AND KEYWORDS
%----------------------------------------------------------------------------------------
%\renewcommand{\abstractname}{Summary} % Uncomment to change the name of the abstract to something else
\begin{abstract}
In the article \cite{Sim}, H. Simmons describes two monads of interests arising from the dual adjunction between the category of topological spaces and that of (bounded) distributive lattices. These are the open prime filter monad and the ideal lattice monad. It is known that the ideal lattice monad induces the ideal frame comonad on the category of frames. We show that this ideal frame comonad can be paired with the open prime filter monad via the open set-spectrum adjunction. From this, we give a new proof of the equivalence between the category of stably compact spaces and that of stably compact frames on one hand, and that of compact Hausdorff spaces and compact regular frames on the other. We show, among other things, how the \v{C}ech-Stone compactification in Pointfree Topology and Pointset Topology relate each other in this particular context.
\end{abstract}

%\hspace*{3,6mm}\textit{Keywords:} lorem , ipsum , dolor , sit amet , lectus % Keywords

\vspace{20pt} % Some vertical space between the abstract and first section

%----------------------------------------------------------------------------------------
%	ESSAY BODY
%----------------------------------------------------------------------------------------

\section{\bf Introduction}

What could be considered as the main results of this paper at first glance are not new, and to a certain extent may be seen as a folklore by some authors. By that, we refer to the equivalence between the category \cat{StKSp} of stably compact spaces and the category \cat{StKFrm} of stably compact frames\footnote{The terminology {\em Stably continuous frames} is also used in the literature (Cf. \cite{BanBru}).} on one hand, and that between the category \cat{KHaus} of compact Hausdorff spaces and the category \cat{KRegFrm} of compact regular frames on the other. The paper attempts to understand these equivalences through a functorial duality between two concepts: filters and ideals. Indeed, on one hand stably compact spaces $X$ appear as the continuous retracts of spaces $FX$ of open  prime filters, and on the other stably compact frames $L$ are retracts, in the category \cat{Frm} of frames and their homomorphisms, of the frame of ideals $\id L$. Beyond the mere fact that prime filters and ideals on the open set lattice $\mathcal{O}X$ complement each other via a lattice homomorphism $\mathcal{O}X\to\set{0,1}$, $F$ and $\id$ respectively underlie a monad \trip{F} on the category \cat{Top} of topological spaces and continuous maps and a comonad \trip{I} on \cat{Frm}. They are paired in the sense that they are part of the commutative diagram
%%%================================================================================================================
$$\xymatrix{
\cat{Top} \ar[d]_{F} \ar[r]^{\mathcal{O}} & \cat{Frm}^{op} \ar[d]^{\id}\\
\cat{Top} & \cat{Frm}^{op} \ar[l]^{\Sigma}
}$$
%%%================================================================================================================
where $\Sigma$ is the spectrum functor. Thus \trip{F} and \trip{I} are dual to each other through the adjunction $\mathcal{O}\dashv\Sigma$. 

Although these two monads are well-known, their studies and use have mainly developed on separate pathways, with a few notable exceptions (Cf. \cite{Ban81}). We note that one such exceptions is H. Simmons' paper ``{\em A Couple of Triples}'' (\cite{Sim}) to which we pay tribute through the title of this paper. The categories \cat{StKSp} and \cat{StKFrm} therefore realize themselves as algebras and coalgebras of these two monads respectively and, with a weaker form of the axiom of choice, their dual equivalence stems from the diagram above. It is a consequence of this same diagram that the categories \cat{KHaus} and \cat{KRegFrm} are dually equivalent. Thus the current work is a contribution to the study of Stone-type dualities.

The categorical language is unavoidable and we have consecrated Section 3 for this purpose. We briefly recall known results on which we build the tools relevant to the construction of the diagram. We consider an adjunction $L\dashv R$, a monad $T$, and compare the two monads $T$ and $RTL$ through their algebras. We show that equivalence is obtained with the assumption $LRT=T$, and that this assumption is equivalent to the Boolean Ultrafilter Theorem (\cite{Ban83}) in the case of the above diagram (Theorem \ref{thm: equivalence of LRT=T with choice}). The main results follow in Section 4 where we first start with Simmons' construction in \cite{Sim}. The last section is devoted to two examples in Topology, as illustrations of the categorical construction and as non-examples for the condition $LRT=T$. The necessary background for a smooth reading and the terminology that we shall use are found in the preliminaries.

\section{\bf Preliminaries}
\paragraph*{\bf Frames}
(\cite{Joh82,Sim}) A distributive lattice  is a partially ordered set\footnote{In the context of the paper, ``set'' means small set.} set admitting finite joins and finite meets, and for which the identity
$$
\meet{a}{(\join{b}{c})}=\join{(\meet{a}{b})}{(\meet{a}{c})}
$$
holds. A homomorphism between distributive lattices is a function that preserves the operations $\wedge$ and $\vee$, including the top element $\wedge\emptyset=1$ and the bottom element $\vee\emptyset=0$. These form the category \cat{DLat}.\\
\noindent
(\cite{PicPul2012,Joh82}) A frame $L$ is a complete lattice with the equational identity
$$
\meet{a}{\left(\ds{\bigvee} S\right)}=\ds{\bigvee\set{\meet{a}{s}\ |\ s\in S}}
$$
for all $a\in L$ and any $S\subseteq L$. Frame homomorphisms or frame maps are functions \map{f}{L}{M} that preserve finite meets  and arbitrary joins. Frames and their homomorphisms form a category denoted by \cat{Frm}. There is a forgetful functor $\cat{Frm}\to\cat{DLat}$. 
\begin{defn}
\label{defn: of an ideal}
A subset $J\subseteq L$ is called an ideal if $J$ is a downset, i.e. $b\in J$ provided there is $a\in J$ with $b\leq a$, and if it is closed under the formation of finite joins. 
\end{defn}
\noindent
The set of all ideals on a distributive lattice $L$ is denoted by $\id L$. $\id L$ is a frame with the following operations: meets are given by set-intersections and joins given by
$$
\bigvee\mathscr{J}=\bigcup\set{I_1\vee I_2\vee\dots\vee I_n\ |\ I_1, I_2, \dots, I_n\in\mathscr{J}\text{ and }n\in \mathbb{N}},
$$
where $I_1\vee I_2\vee\dots\vee I_n=\set{i_1\vee i_2\vee\dots\vee i_n\ |\ i_k\in I_k\text{ for }1\leq k\leq n}$. Thus, if $\mathscr{J}$ is a directed set, then $\ds{\bigvee}\mathscr{J}=\bigcup\mathscr{J}$. If \map{f}{L}{M} is a  homomorphism between distributive lattices, then the function \map{\id f}{\id L}{\id M} given by
$$
\id f(I)=\set{b\ |\ b\leq f(a)\text{ for some }a\in I},
$$
is a frame homomorphism. In particular \id\ can be viewed as an endofunctor on \cat{DLat} with range in \cat{Frm}. The frames of the form $\id D$, where $D$ is a distributive lattice, are called {\em coherent frames} (\cite{Ban81,Joh82}). \\
\noindent
Consider the category $\cat{Loc}=\cat{Frm^{op}}$ of locales or frames. Extracting the topology $\mathcal{O}X$ of a space $X$, and taking $\mathcal{O}(f)=\inv{f}[-]$ for each continuous function $f$, gives a functor \map{\mathcal{O}}{\cat{Top}}{\cat{Loc}}. $\mathcal{O}$ admits a right adjoint \map{\Sigma}{\cat{Loc}}{\cat{Top}} which is defined as follows: for each frame $L$, one can assign the space $\Sigma L=\set{\map{f}{L}{2}\ |\ f\text{ is a frame map}}$ whose basic open sets are given by \set{\Sigma_a\ |\ a\in L}, where $\Sigma_a=\set{\map{f}{L}{2}\ |\ f(a)=1}$ consists of all points contained in $a\in L$. The adjunction $\mathcal{O}\dashv\Sigma$ reduces itself to a categorical equivalence between sober spaces and partial frames.\\

\noindent
\paragraph*{\bf Stably compact frames} The way below relation $\ll$ on a frame $L$ is defined as follows: $a\ll b$ if for any $S$ with $b\leq \bigvee S$, there is a finite $G\subseteq S$, such that $a\leq \bigvee G$. For open sets $U$ and $V$ in a space, the relation $U\ll V$ means that $U$ is relatively compact in $V$. A frame $L$ is sadi to be {\em stably compact} if $\ll$ is a sublattice of $L\times L$ and for any $a\in L$:
\begin{center}
$a=\bigvee\set{x\ |\ x\ll a}$.
\end{center}
{\em Stably continuous frame} is another name for stably compact frame. Here we emphasize the fact that $1\ll 1$.  However, the two are sometimes interchangeably used in the literature (Cf. \cite{BanBru}). A {\em stably compact space} is compact, locally compact, sober and stable - i.e. the relation $\ll$ on $\mathcal{O}X$ is stable. There are various equivalent definitions of a stably compact space (See \cite[References]{JunKegMos2001}). However for the sake of brevity, we choose the description that is close to the pointfree definition.\\

\paragraph*{\bf Natural transformations.}The composition of morphisms $f$ and $g$ will be denoted by \comp{g}{f} or just simply by $gf$. Given two pairs of functors \map{F,G}{\cat{C}}{\cat{D}} and \map{H,K}{\cat{D}}{\cat{E}} and two natural transformations $\alpha:F\to G$ and $\beta:H\to K$, the horizontal composition $\beta\circ\alpha:HF\to KG$ is defined by 
\begin{center}
$(\beta\circ\alpha)_X=\beta_{GX}H(\alpha_X)=K(\alpha_X)\beta_{FX}$
\end{center}
for each $X$ in \cat{C}. This is usually written as $\beta\circ\alpha=\beta G\cdot H\alpha=K\alpha\cdot\beta F$.
If \map{L}{\cat{C}}{\cat{D}} is another functor and $\delta:G\to L$ another natural transformation, then the vertical composition $\comp{\delta}{\alpha}:F\to L$ is simply defined as $\comp{\delta_X}{\alpha}_X$ for each $X\in\cat{C}$. The two types of compositions interact as follows  
\begin{center}
$(\beta\smash'\circ\alpha\smash')\cdot(\beta\circ\alpha)=(\beta\smash'\cdot\beta)\circ(\alpha\smash'\cdot\alpha)$. (See \cite{Low,Mac}.)
\end{center}
\noindent
For a category \cat{C} and two objects $A$ and $B$, we denote by $\cat{C}(A,B)$ the set of morphisms from $A$ to $B$.\\

\paragraph*{\bf Monads.} Recall that a {\em monad} $\mathbb{T}$ on a category \cat{C} is a triple \mon{T}{m}{e}, where \map{m}{TT}{T} and \map{e}{1}{T} are natural transformations satisfying the identities 
\begin{center}
$\comp{m}{Tm}=\comp{m}{m T}$ and $\comp{m}{e T}=\comp{m}{Te}=1_T$.
\end{center}
A \trip{T}-algebra (or an {\em Eilenberg-Moore algebra}) is a pair $(X,a)$, where $X\in\cat{C}$ and \map{a}{TX}{X} a morphism such that 
\begin{center}
$\comp{a}{Ta}=\comp{a}{m_X}$ and $\comp{a}{e_X}=1_X$.
\end{center}
A morphism between two monads $\trip{T}=\mon{T}{m}{e}$ and $\trip{M}=\mon{M}{n}{d}$ is a natural transformation \map{\alpha}{T}{M} satisfying $\comp{\alpha}{e}=d$ and $\comp{\alpha}{m}=\comp{n}({\ocomp{\alpha}{\alpha}})$. 

\noindent 
Note that $(TX,m_X)$ is the free \trip{T}-algebra over $X$. If $(X,a)$ and $(Y,b)$ are \trip{T}-algebras, then a \trip{T}-algebra homomorphism \map{f}{(X,a)}{(Y,b)} is a morphism \map{f}{X}{Y} in \cat{C} such that $\comp{f}{a}=b\cdot Tf$. The category of \trip{T}-algebras and \trip{T}-algebra homomorphisms are denoted by \powcat{\cat{C}}{T} or by \alg{T}. The forgetful functor $\map{\powcat{G}{T}}{\powcat{\cat{C}}{T}}{\cat{C}}:(X,a)\mapsto X$ admits a left adjoint $\map{\powcat{F}{T}}{\cat{C}}{\powcat{\cat{C}}{T}}:X\mapsto (TX,m_X), f\mapsto Tf$. The unit of this adjunction is given by $e_X:X\to \powcat{G}{T}\powcat{F}{T}X$ and the co-unit is provided by the algebra morphisms $\varepsilon_{(X,a)}:\powcat{F}{T}\powcat{G}{T}(X,a)\to(X,a)$. 

\noindent 
Any adjunction $(F\dashv G,\eta,\epsilon)$ where \map{G}{\cat{C}}{\cat{A}} and \map{F}{\cat{A}}{\cat{C}} induces a monad $(GF,G\epsilon F,\eta)$ on \cat{C}. If the adjunction $F\dashv G$ induces the monad \trip{T} on \cat{C}, that is $GF=T, G\epsilon F=m$ and $\eta=e$, then there is a unique functor\map{K}{\cat{A}}{\powcat{\cat{C}}{T}} - called {\em comparison functor}, such that $KF=\powcat{F}{T}$ and $\powcat{G}{T}K=G$:
%%%%====================================================================
$$\xymatrix{
& & \powcat{\cat{C}}{T} \ar@<1.0ex>[dd]^{\small{\powcat{G}{T}}}_{\dashv}\\
 & & \\
\cat{A} \ar@{-->}[uurr]^{K} \ar@<-1.0ex>[rr]_{G} &  & \cat{C} \ar@<-1.1ex>[ll]_{F}^{\perp} \ar@<1.0ex>[uu]^{\powcat{F}{T}} 
}$$
%%%%====================================================================
\noindent
$G$ is said to be monadic if $K$ is an equivalence. In that case, the relations $G(X,\alpha)\cong G(Y,\beta)$ and $(X,\alpha)\cong (Y,\beta)$ are equivalent.

\section{\bf Lifting a monad through an adjunction}
We consider a monad $\trip{T}=(T,m,e)$ on \cat{C} and an adjunction $(L\dashv R,\eta,\epsilon)$ where \map{R}{\cat{C}}{\cat{B}} and \map{L}{\cat{B}}{\cat{C}}. Since $\powcat{F}{T}L\dashv R\powcat{F}{T}$ the composition $RTL=R\powcat{G}{T}\powcat{F}{T}L$ underlies a monad $\trip{M}=(M,n,d)$ on \cat{B} with multiplication $n=RmL\cdot RT\epsilon TL$ and unit $d=\comp{ReL}{\eta}$. 
We start with a simple example.
\begin{ex}
\label{ex: closure operators}
Let \map{f}{X}{Y} be a continuous function and let us denote by $c_X$ and $c_Y$ the Kuratowski closure operators on $X$ and $Y$ respectively. Let $Cl(X)=\set{A\ |\ c_X(A)=A}$ and let $Cl(Y)=\set{B\ |\ c_X(B)=B}$. $f$ induces an image-pre-image adjunction $f^*\dashv f_*$ on the subobjects of $X$ and $Y$. Considering the lattices of subobjects as categories and the closure operators as monads, the composition $f^*c_Yf_*$ gives another closure operator on $X$. When $c_X=f^*c_Yf_*$, we say that $f$ is $c$-{\em initial}. The use of initial morphisms, along final morphisms, finds its importance in the study of topological connectedness in categories (\cite{Cle2001}). 
\end{ex}
\noindent
Now, consider the comparison functor \map{K}{\powcat{\cat{C}}{T}}{\powcat{\cat{B}}{M}}:
%%%%====================================================================
$$\xymatrix{
\powcat{\cat{C}}{T} \ar@{-->}[rr]^K \ar@<1.0ex>[dd]^{\small{\powcat{G}{T}}}_{\dashv} & & \powcat{\cat{B}}{M} \ar@<1.0ex>[dd]^{\small{\powcat{G}{M}}}_{\dashv}\\
 & & \\
\cat{C} \ar@<1.0ex>[uu]^{\small{\powcat{F}{T}}}  \ar@<-1.0ex>[rr]_{R} &  & \cat{B} \ar@<-1.1ex>[ll]_{L}^{\perp} \ar@<1.0ex>[uu]^{\small{\powcat{F}{M}}} 
}$$
%%%%====================================================================
Since $G\powcat{G}{T}=\powcat{G}{M}K$, the functor $K$ is an {\em algebra lifting} of $R$, that is, $R$ restricts to a functor between the algebras. We recall some classical results about such a lifting even when the left adjoint $L$ of $R$ is absent.
\begin{thm}
\label{thm: adjunction lifting}
(\cite{App,Bor2,Joh75,ManMul}) There is natural bijection between liftings \map{K}{\powcat{\cat{C}}{T}}{\powcat{\cat{B}}{M}} and natural transformation \map{\lambda}{MR}{RT} satisfying
%%%%====================================================================
$$\xymatrix{
R\ar[r]^{dR} \ar[dr]_{Re} & MR \ar[d]^{\lambda} & MMR \ar[l]_{nR} \ar[d]^{M\lambda} \\
 & RT & MRT \ar[d]^{\lambda T} \\
 & & RTT \ar[lu]^{Rm}
}$$
%%%%====================================================================
\end{thm}
%\begin{thm}
%\label{thm: existence of left adjoint for lifting}
%(\cite{Bor2,Joh75}) If \powcat{\cat{C}}{T} has coequalizers of reflexive pairs and the left adjoint $L$ exists, then $K$ admits a left adjoint.
%\end{thm}
\noindent 
Since \trip{M} is obtained from $L,R$ and $T$, these statements can be simplified further and an explicit expression for the law $\lambda$ can be obtained.
\begin{prop}
\label{prop: comparison functor as a lifting}
The natural transformation \map{RT\eps}{(RTL)R}{RT} satisfies the diagram in Theorem \ref{thm: adjunction lifting}.
\end{prop}
\begin{proof}
Naturality of \eps\ gives
%%%%====================================================================
$$\xymatrix{
RTLR \ar[dd]_{\ds{RT\eps}} & & & & RTLRTLR \ar[llll]_{\left(\ds{\comp{Rm}{RT\eps T}}\right)\ds{LR}} \ar[d]^{\ds{RTLRT\eps}} \\
 & & & & RTLRT \ar[d]^{\ds{RT\eps T}} \\
RT & & & & RTT \ar[llll]^{Rm}
}$$
%%%%====================================================================
The other diagram is given by
%%%%====================================================================
$$\xymatrix{
R\ar[r]^{\eta R} \ar[dr]_{1} & RLR \ar[d]^{\ds{R\eps}} \ar[rr]^{ReLR} & & RTLR \ar[d]^{\ds{RT\eps}}\\
& R \ar[rr]_{Re} & & RT 
}$$
%%%%====================================================================
where $\comp{ReLR}{\eta}R=dR$.
\end{proof}
\begin{prop}
\label{prop: comparison functor as a mapping}
The comparison functor $K$ corresponds to the natural transformation law $RT\eps$.
\end{prop}
\begin{proof}
As a lifting, $K$ assigns to a \trip{T}-algebra $(X,\alpha)$ the composition \map{\comp{R\alpha}{\lambda_X}}{RTLRX}{RX}, where \map{\lambda}{RTLR}{RT} is the unique law associated to $K$. Since $K$ is a comparison functor, $\powcat{G}{M}K$ takes $(X,\alpha)$ to the $R$-image of the counit 
%%%%====================================================================
$$\xymatrix{
 \powcat{F}{T}LR\powcat{G}{T}(X,\alpha)\ar[rrr]^{\powcat{F}{T}\eps\powcat{G}{T}}& & & \powcat{F}{T}\powcat{G}{T}(X,\alpha)\ar[r]^{\ds{\alpha}} & (X,\alpha)
}$$
%%%%====================================================================
of the adjunction $\powcat{F}{T}L\dashv R\powcat{G}{T}$. This is explicitly given by
%%%%====================================================================
$$\xymatrix{
 (TLRX,m_{LRX})\ar[rrr]^{\ds{T(\eps_X)}} & & & (TX,m_X)\ar[r]^{\ds{\alpha}} & (X,\alpha)
}$$
%%%%====================================================================
We have $\comp{R\alpha}{\lambda_X}=\comp{R\alpha}{R\ds{T(\eps_X)}}$. Therefore $\lambda=RT\eps$.
\end{proof}
\noindent 
We give further useful observations. Let $\cat{Mon(D)}$ denote the category of monads over a category \cat{D} together with morphisms of monads.
\begin{prop}
\label{prop: lifting of monads is functorial}
The assignment \map{\mathfrak{L}}{\cat{Mon(C)}}{\cat{Mon(B)}}, defined by $\mathfrak{L}(T)=RTL$ on objects and by $\mathfrak{L}(\map{\lambda}{T}{N})=R\lambda L$ on morphisms, is well-defined and functorial. 
\end{prop}

\begin{proof}
Functoriality of \lift{L} is straightforward. Let \map{\lambda}{(T,m,e)}{(N,n,d)} be a morphism of monads. Since $\comp{\lambda}{e}=d$, we have $\comp{R\lambda L}{(R eL\cdot\eta)}=R d L\cdot\eta$. In the following diagram
%%%%%=================================================================================================
$$\xymatrix{
TLRT\ar[rr]^{\ds{T\eps T}} \ar[d]_{\ds{TLR\lambda}} & & TT \ar[d]^{T\lambda} \ar[r]^m & T\ar[dd]^{\lambda}\\
TLRN \ar[rr]^{\ds{T\eps N}} \ar[d]_{\ds{\lambda LRN}} & & TN \ar[d]^{\lambda N} & \\
NLRN \ar[rr]_{\ds{N\eps N}} & & NN\ar[r]_{n} & N
}$$
%%%%%=================================================================================================
where $\ocomp{\lambda}{\lambda}=\comp{\lambda N}{T\lambda}$, we have $\lambda\cdot m\cdot T\eps T=n\cdot N\eps N\cdot \lambda LRN\cdot TLR\lambda$. Now, composing with $R\cdot (-)\cdot L$ gives 
\begin{center}
$\comp{R\lambda L}{(\comp{RmL}{RT\eps TL})}=\comp{\comp{RnL}{RN\eps NL}}{(\ocomp{R\lambda L}{R\lambda L})}$. 
\end{center}
Thus $R\lambda L$ is a morphism of monads. 
\end{proof}

\begin{cor}
\label{cor: monad morphism between the two induced monads}
The natural transformation \map{ReL}{RL}{RTL} is a morphism of monads.
\end{cor}
\begin{proof}
The natural transformation \map{e}{1}{T} is a morphism of monads.
\end{proof}
\begin{prop}
\label{prop: functor between two induced algebras and RTL-algebras are RL-algebras}
Denote $\trip{H}=(RL,R\eps L,\eta)$. The morphism of monads $ReL$ corresponds to an embedding functor \map{V}{\powcat{\cat{B}}{M}}{\powcat{\cat{B}}{H}} such that $\powcat{G}{H}V=\powcat{G}{M}$. 
\end{prop}
\noindent
It should not be expected that $K$ is generally an equivalence. In Example \ref{ex: closure operators}, the concrete case of the inclusion function \map{f}{(0;1)}{\R} shows that the restriction of the pre-image function $f^*$ to a monotone function $Cl(\R)\to Cl(0;1)$ is not faithful. It is however known (See \cite{Joh75} and \cite{Lin}) that if \powcat{\cat{C}}{T} has coequalizers of reflexive pairs, then $K$ admits a left adjoint. A slightly stronger but similar observation holds for the functor \map{V}{\powcat{\cat{B}}{M}}{\powcat{\cat{B}}{H}}.
\begin{thm}
\label{thm: reflexive coequalizers imply existence of left adjoint}
(\cite{Joh75,Lin}) If \powcat{\cat{C}}{T} has coequalizers of reflexive pairs, then $K$ admits a left adjoint. 
\end{thm}
\noindent
Similarly, we have
\begin{prop}
\label{prop: V monadic}
(\cite[Corollary 4.5.7]{Bor2}) With the monads \trip{H} and \trip{M}, and functor $V$ as in Proposition \ref{prop: functor between two induced algebras and RTL-algebras are RL-algebras}, if the category \powcat{B}{M} has coequalizers, then $V$ is monadic.
\end{prop}
\noindent 
The general context in which Theorem \ref{thm: reflexive coequalizers imply existence of left adjoint} above appears in the papers \cite{Joh75} and \cite{Lin} does not assume a prior existence of the left adjoint $L$. This is only added in the statement of the results. According to Johnstone, the first explicit appearance of this result was in the Appendix to Diaconescu's thesis \cite{Dia}. When \cat{C} admits a proper $(\mathcal{E},\mathcal{M})$-factorisation system, mild conditions on \cat{C} and $T$ are sufficient for \powcat{\cat{C}}{T} to be cocomplete.

\begin{thm}
\label{thm: cocompleteness of algebras with (E,M) fact system}
(\cite{Ada}) Suppose that \cat{C} is cocomplete and $\mathcal{E}$-well-copowered. If $T$ preserves $\mathcal{E}$, that is $Th\in\mathcal{E}$ for all $h\in\mathcal{E}$, then \powcat{\cat{C}}{T} is cocomplete.
\end{thm}
\noindent 
Weaker versions of this theorem also appear in \cite{Lin} and \cite{Bar70b} where completeness of \cat{C} is also considered. In addition, as noted in \cite{Ada}, Linton assumed that $T$ preserves the class $\mathcal{M}$ as well. Preservation of the class $\mathcal{E}$ seems to be essential for monads and reflectors that arise from the study of compactifications. (Cf. \cite{Hol00} and \cite{Raz}.)\footnote{Reflectors are essentially idempotent monads. (Cf. \cite{Bor2}.) Note that in \cite{Hol00}, preservation of $\mathcal{E}$ means that a reflector sends a morphism in $\mathcal{E}$ into a morphism in $\mathcal{E'}$, where $\mathcal{E'}$ is part of the factorisation system in the subcategory of concern.} We shall show that under special circumstances, it is possible to show that \powcat{\cat{B}}{M} is equivalent to a reflective subcategory of \powcat{\cat{C}}{T}.
\begin{lem}
\label{lem: Every T-algebra is a fixed point of LR and eta is an iso}
If $LRT\cong T$, then $X\cong LRX$ for any given \trip{T}-algebra $(X,\alpha)$. In this case, there is an embedding \map{E}{\powcat{\cat{C}}{T}}{\cat{Fix}(LR)} that is monadic.
\end{lem}
\begin{proof}
First note that if  $LRT\cong T$, then also $\eps T$ is a natural isomorphism by virtue of its universality. It is straightforward to verify that the inverse of \map{\eps_X}{LRX}{X} is given by $LR\alpha\cdot\inv{\eps_{TX}}\cdot \eta_X$.
\end{proof}
\begin{lem}
\label{lem: eta is iso for RTL algebra if LRT=T}
Suppose that $LRT\cong T$, then for any $RTL$-algebra $(X,\alpha)$ the morphism \map{\eta}{X}{RLX} is an isomorphism. This gives an embedding \map{E'}{\powcat{\cat{B}}{M}}{\cat{Fix}(RL)} that is monadic.
\end{lem}
\begin{proof}
We have $R\eps TL\cdot \eta RTL=1_{RTL}$ and since $R\eps TL$ is an isomorphism, so is $\eta RTL$ with inverse $\inv{(\eta RTL)}=R\eps TL$. Thus, from the natural square $RL\alpha\cdot\eta RTL=\eta\cdot\alpha$, we have $1=\eta\cdot\alpha\cdot R\eps TL\cdot RL(ReL\cdot\eta)$. Therefore $\eta$ is an isomorphism.
\end{proof}
\begin{thm}
\label{thm: existence of a left adjoint when RLT=T}
If  $LRT\cong T$, then the comparison functor $K$ is an equivalence.
\end{thm}
\begin{proof}
 Define a functor \map{K^*}{\powcat{\cat{B}}{M}}{\powcat{\cat{C}}{T}} as follows: for each \trip{M}-algebra $(Y,\beta)$, take $K^*(Y,\beta)=(LY,\comp{L\beta}{\inv{(\eps TL_Y)}})$ with the obvious assignment on morphisms. The right transpose of $eL_Y$ is $\comp{ReL_Y}{\eta_Y}$ and so we have $\comp{\eps TL_Y}{L(\comp{ReL_Y}{\eta_Y})}=eL_Y$. We have
 \begin{align*}
 (\comp{L\beta}{\inv{(\eps TL_Y)}})\cdot eL_Y &=\comp{L\beta}{L(\comp{ReL_Y}{\eta_Y})}\\
 &=L(\comp{\beta}{\comp{ReL_Y}{\eta_Y}})\\
 &=L(1_Y)=1_{LY}.
 \end{align*}
The following diagram 
%%%===================================================================================================
$$\xymatrix{
TTLY \ar[d]_{\inv{(T\eps TL)}} \ar[rrrrr]^{\comp{TL\beta}{\inv{(T\eps TL)}}} & & & & &  TLY \ar[d]_{\inv{(\eps TL)}}\\
TLRTLY\ar[d]_{\ds{T\eps TL}} & & & LRTLRTLY\ar[lll]_{\ds{\eps TLRTL}} \ar[d]^{\ds{LRT\eps TL}}\ar[rr]^{LRTL\beta} & &  LRTLY\ar[dd]^{L\beta} \\
TTLY \ar[d]_{mL} & & & LRTTLY\ar[lll]_{\ds{\eps TTL}}\ar[d]^{LRmL} &  & \\ 
TLY & & & LRTLY\ar[lll]^{\ds{\eps TL}}\ar[rr]_{L\beta} & &  LY
}$$
%%%===================================================================================================
where $n=RmL\cdot RT\eps L$, shows that $L\beta\cdot\inv{(\eps TL)}$ is indeed a \trip{T}-algebra morphism. Now, for any $RTL$-algebra $(Y,\beta)$ we have $KK^*(Y,\beta)=(RLY,\comp{RL\beta}{\comp{R(\inv{\eps TL})}{RT\eps L}})$. By Beck's monadicity theorem and by Lemma \ref{lem: eta is iso for RTL algebra if LRT=T}, $\eta$ is an isomorphism of $RTL$-algebras and so $KK^*\cong 1$. Similarly, by Lemma \ref{lem: Every T-algebra is a fixed point of LR and eta is an iso}, $\eps$ is an isomorphism of $T$-algebras.
\end{proof}
\noindent
One can easily check that with the assumptions of the above result, $K^*$ satisfies $\powcat{F}{T}L\cong K^*\powcat{F}{M}$. The assumption $LRT\cong T$ has some further implication.
\begin{prop}
\label{prop: monoidal composition of lifting}
Let $T$ be a monad such that $LRT\cong T$. If $NT$ is a monad, then $\lift{L}(NT)=\lift{L}(N)\lift{L}(T)$. In particular, if such a $T$ is a reflector, i.e. $T$ is idempotent, then $\lift{L}(T)$ is also a reflector.
\end{prop}

\section{\bf The open prime filter monad and the ideal frame comonad}
Taking the lattices of ideals $\id L$ induces a comonad on \cat{Frm}, hence a monad with an underlying functor $\overline{\id}$ on \cat{Loc}. This gives a monad with underlying functor $\Sigma\overline{\id}\mathcal{O}$ on \cat{Top}. Although this monad is described by Simmons in \cite{Sim} and separately by Br\"ummer in \cite{Bru79}, we shall briefly discuss how it was generated as we need to show through the results shown by Banaschewski in \cite{Ban81} how it relates to the functor $\id$ on frames.

\noindent 
Consider the  Sierpinski space $(\set{0,1},\set{\emptyset,\set{1},\set{0,1}})$ and the distributive lattice $\set{0< 1}$. Both objects are manifestations of the schizophrenic object $2$ with which we denote them. The adjunction (\cite{Sim})
\begin{center}
\Adj{\cat{Top}}{\cat{Top}(-,2)}{\cat{DLat}(-,2)}{\cat{DLat}^{op}}
\end{center}
induces a monad $\mathbb{F}=(F,\mu,\eta)$ on \cat{Top}, where $FX$ is the set of open prime filters on $(X,\tau)$, whose topology admits the basis \set{O^{*}\ |\ O\in\tau}, and where
\begin{center}
$\mu(\lift{X})=\set{O\in\tau\ |\ O^{*}\in\mathfrak{X}}$ and $\eta(x)=\set{O\in\tau\ |\ x\in O}$.
\end{center}
On the other hand, the ideal lattice monad $\mathbb{I}=(\id,\bigcup,\downarrow)$ is induced by the same ajdunction on \cat{DLat}, where
\begin{center}
$\map{\bigcup}{\id\id D}{\id D}: \mathscr{I}\mapsto \bigcup\mathscr{I}$ and $\downarrow(a)=\set{x\ |\ x\leq a}$.
\end{center}
for each distributive lattice $D$.
\begin{thm}
\label{thm: algebras of F and I}
\begin{enumerate}
\item (Cf. \cite[Section 3]{Sim} and \cite[Proposition 5]{Ban81}) The algebras of the monad $\mathbb{F}$ are precisely the stably compact spaces and proper maps.
\item (Cf. \cite[Section 2]{Sim} and \cite[Exercise VI.4.6]{Joh82}) The algebras of the monad $\mathbb{I}$ are precisely the frames and frame homomorphisms.
\end{enumerate}
\end{thm}
\noindent
The adjunction between \cat{DLat} and its $\mathbb{I}$-algebras \cat{Frm} in turn induces a comonad $\mathbb{K}=(\id,c,\gamma)$ on \cat{Frm} where $c=\id(\downarrow)$ and $\gamma$ are defined as follows:
\begin{center}
$\map{c}{\id L}{\id\id L}: I\mapsto \set{J\ |\ \bigvee J\in I}$ and $\gamma(a)=\set{x\ |\ x\ll a}.$
\end{center}
\begin{thm}
\label{thm: co algebras of K}
(Cf. \cite[Section 3]{BanBru}) The coalgebras of the comonad $\mathbb{K}$ are exactly the stably compact frames and proper frame homomorphisms.
\end{thm}
\noindent
From \cite[Lemma 2]{Ban81}, the natural equivalence $\cat{Frm}(\id D,2)\cong\cat{DLat}(D,2)$ induces a natural equivalence between $\Sigma\cdot \id$ and $\cat{DLat}(-,2)$. This is because for a distributive lattice $D$, the completely prime filters of $\id D$ are in natural bijection with prime filters of $D$. Thus we have for a space $X$ that $(\Sigma\cdot\id\cdot\mathcal{O})(X)=FX$. To show that this indeed gives the same monad $\mathbb{F}$, consider the following diagram:
%%%================================================================================================================
$$\xymatrix{
 & & \cat{DLat}^{op}\ar@<0.55ex>[dll]^{\pi} \ar@<1ex>[dd]^{F^\mathfrak{I}} \\
\cat{Top}\ar@<1.75ex>[urr]^{\Omega}_*-<0.55ex>{_\perp} \ar@<1.25ex>[drr]^{\mathcal{O}}_*-<0.55ex>{_\perp} & & \\
 & & \cat{Frm}^{op} \ar@<1ex>[uu]^{G^{\mathfrak{I}}}_*-<0.55ex>{_\dashv} \ar@<1ex>[ull]^{\Sigma}
}$$
%%%================================================================================================================
where $\pi=\cat{DLat}(-,2)$ and $\Omega=\cat{Top}(-,2)$. For the vertical arrows, $\comp{G^{\mathfrak{I}}}{F^{\mathfrak{I}}}=\comp{\Omega}{\pi}$ is the lattice ideal monad on \cat{DLat} and $\comp{F^{\mathfrak{I}}}{G^{\mathfrak{I}}}=\id$ is the ideal frame comonad on \cat{Frm}. The functor $\mathcal{O}$ acts as comparison functor between $\Omega\dashv\pi$ and $G^{\mathfrak{I}}\dashv F^{\mathfrak{I}}$ so that $\comp{G^{\mathfrak{I}}}{\mathcal{O}}=\Omega$ and $\comp{\mathcal{O}}{\pi}=F^{\mathfrak{I}}$. 
\begin{prop}
\label{prop: open prime filter monad innduced from ideal frame comonad}
The monad $\mathbb{F}$ is induced by the adjunction $\mathcal{O}\dashv\Sigma$ from the comonad $\mathbb{K}=(\id,c,\gamma)$.
\end{prop}
\begin{proof}
We have
\begin{center}
$\Sigma\cdot \id\cdot \mathcal{O} = \Sigma\cdot (\comp{F^{\mathfrak{I}}}{G^{\mathfrak{I}}})\cdot \mathcal{O}$,
\end{center}
and so $\Sigma\cdot\id\cdot\mathcal{O}$ is induced by the adjunction $G^{\mathfrak{I}}\cdot \mathcal{O}\dashv \Sigma\cdot F^{\mathfrak{I}}$. Since right adjoints are unique up to equivalence\footnote{This actually manifests itself in the natural equivalence $\cat{Frm}(\id D,2)\cong\cat{DLat}(D,2)$.}, we have $\pi\cong \comp{\Sigma}{F^{\mathfrak{I}}}$, and so $\Sigma\cdot \id\cdot \mathcal{O}\cong\comp{\pi}{\Omega}$. Therefore $\Sigma\cdot \id\cdot \mathcal{O}$ induces the same open prime filter monad $\mathbb{F}$ on \cat{Top}.
\end{proof}
\noindent
We will now proceed to describe the relationship between stably compact spaces and stably compact frames in this context.
\begin{prop}
\label{prop: spectrum restriction admits a left adjoint}
The category $\cat{StkFrm}$ of stably compact frames is complete and the category \cat{StKSp} of stably compact spaces is essentially coreflective inside the category $\cat{StkFrm}^{op}$. 
\end{prop}
\begin{proof}
First, the category \cat{Frm} is complete. It has an $(\mathcal{E},\mathcal{M})$-factorisation system where $\mathcal{E}$ is the class of extremal (or regular) epimorphisms and $\mathcal{M}$ the class of monomorphisms which are injections (\cite[Section III.1.4]{PicPul2012}). \cat{Frm} is $\mathcal{M}$-wellpowered, and the functor \id\ preserves the class $\mathcal{M}$. By Theorem \ref{thm: reflexive coequalizers imply existence of left adjoint} and Theorem \ref{thm: cocompleteness of algebras with (E,M) fact system}, $\cat{StkFrm}$ is complete and the restriction \map{K}{\cat{StKFrm}^{op}}{\cat{StKSp}} of $\Sigma$ admits a left adjoint $K^*$ which - as a restriction of $\mathcal{O}$ - is necessarily a full embedding.
\end{proof}
\noindent
When the equation $LRT\cong T$ in Theorem \ref{thm: existence of a left adjoint when RLT=T} is involved, we reconnect with the results of Banaschewski in \cite{Ban83}. In particular:
\begin{thm}
\label{thm: equivalence of LRT=T with choice}
(\cite[Proposition 0]{Ban83}) The following are equivalent:
\begin{enumerate}
\item Every coherent frame is spatial;
\item The Boolean Ultrafilter Theorem; and
\item $\mathcal{O}\cdot \Sigma\cdot\id\cong\id$.
\end{enumerate}
\end{thm}
\begin{proof}
The equivalence between (1) and (3) follows from the observation that a coherent frame is already stably compact and hence belongs to the fixed objects of the coreflection in Proposition \ref{prop: spectrum restriction admits a left adjoint} once it is spatial.
\end{proof}
\begin{cor}
\label{cor: equivalence between stksp and stkfrm}
If the Boolean Ultrafilter Theorem holds, then the categories \cat{StKSp} and \cat{StKFrm} are dually equivalent.
\end{cor}
\noindent
We mention the following result which is also known.
\begin{prop}
\label{prop: monadicity if subcategories}
With the Boolean Ultrafilter Theorem, the category of \cat{StKSp} is monadic inside the category of sober spaces and continuous maps, and \cat{StKFrm} is comonadic inside the category spatial frames and frame homomorphisms.
\end{prop}
\noindent
The restriction of the open prime filter monad $\mathbb{F}$ provides a maximal $T_0$ stable compactification of sober spaces (See \cite{BezHar2014,BezHar2023,Raz,Sal00}), analogous to role played by the \v{C}ech-Stone compactification for Tychonoff spaces. On the other hand, because of Proposition \ref{prop: spectrum restriction admits a left adjoint}, \cat{StKSp} is cocomplete and so, by Proposition \ref{prop: V monadic}, is also monadic inside the category of $T_0$ spaces and continuous maps. This is not very surprising since $\mathbb{F}$ is usually defined on the later (\cite{Hof2013}).\\
\noindent
The equivalence in Corollary \ref{cor: equivalence between stksp and stkfrm} reduces to an equivalence between the category \cat{KHaus} of compact Hausdorff spaces and continuous maps and the (dual) category \cat{KRegFrm} of compact regular frames and frame homomorphisms. Consider the following diagram where only the lower rectangle does not commute:
%%%=================================================================================================================
$$\xymatrix{
\cat{Top}\ar[rr]^{\mathcal{O}}\ar[d]_{F} & & \cat{Frm}^{op}\ar[d]^{\id}\\
\cat{Top}\ar@/_1.0pc/@{-->}[rr]_{\mathcal{O}}\ar[d]_{R} & & \cat{Frm}^{op}\ar@/_1.0pc/@{-->}[ll]_{\Sigma}\ar[d]^{CReg}\\
\cat{Top} & & \cat{Frm}^{op}\ar[ll]^{\Sigma}
}$$
%%%=================================================================================================================
Here $R$ is the Hausdorff reflector, and $CReg$ is the coreflector that assigns to a given frame $L$ is coreflector $CReg L$ (\cite{Joh82,PicPul2012}). It is known that the composition $R\cdot F$ is the \v{C}ech-Stone compactification reflector (\cite{Raz,Sal00}) and that the composition $CReg\cdot\id$ is also the pointfree \v{C}ech-Stone compactification coreflector (\cite{Joh82,PicPul2012}) on \cat{Frm}.
\begin{prop}
\label{prop: KHaus coreflective in KRegFrm}
\cat{KRegFrm} is complete and \cat{Khaus} is essentially coreflective in the category $\cat{KRegFrm}^{op}$.
\end{prop}
\begin{rem}
\label{rem: cocompleteness of KHaus}
Completeness of \cat{KRegFrm} and \cat{StKFrm}, as well as cocompleteness of \cat{Khaus} and \cat{StKSp}, may seem counter-intuitive. To understand this, consider the countable family $\set{1},\set{2}$, $\dots$, $\set{n}$, $\dots$ in \cat{KHaus}. Their sum is given by $\mathbb{N}$, but their colimit in \cat{KHaus} is provided by the \v{C}ech-Stone compactification $\beta\N$.
\end{rem}
\begin{prop}
\label{prop: equivalence between compact Hausdorff spaces and compact regular frames}
If the Boolean Ultrafilter Theorem holds, then the dual equivalence in Corollary \ref{cor: equivalence between stksp and stkfrm} reduces to a dual equivalence between \cat{KHaus} and \cat{KRegFrm} as algebras of $R\cdot F$ and coalgebras of $CReg\cdot\id$ respectively.
\end{prop}
\begin{proof}
Since each $FX$ is sober - hence $T_0$, $\Sigma\cdot Creg\cdot\mathcal{O}\cdot F=R\cdot F$. By Proposition \ref{prop: monoidal composition of lifting}, $\lift{L}(CReg\cdot\id)=\lift{L}(CReg)\cdot\lift{L}(\id)\cong R\cdot F$. By the Boolean Ultrafilter Theorem, $\mathcal{O}\cdot\Sigma$ fixes the range of $CReg\cdot\id$ (See \cite[Proposition 0]{Ban83} and \cite[Proposition 5 and Remark]{Ban81}). The result then follows from Theorem \ref{thm: equivalence of LRT=T with choice}.
\end{proof}

\begin{rem}
\label{rem: about T_0 and commutative lower triangle}
For a $T_0$ space $X$, $X$ is compact and Hausdorff if and only if $\mathcal{O}X$ is compact and regular. Thus the equality $\comp{\Sigma}{\comp{CReg}{\mathcal{O}}}=R$ does not necessarily holds on its own but only with respect to the composition with the open prime filter functor $F$, for which each $FX$ is sober, hence $T_0$. If instead we consider the ultrafilter monad $\mathbb{U}$ with underlying functor \lift{U} on \cat{Top} (\cite[Section 5.1]{Raz}), we still have $\lift{L}(CReg\cdot\id)\cong\comp{R}{\mathfrak{U}}=\comp{R}{F}$. However, although $\comp{\mathcal{O}}{\mathfrak{U}}=\comp{\id}{\mathcal{O}}$, one has $\lift{L}(\id)=F\neq\mathfrak{U}$. In fact, the monad $\mathbb{F}$ is the $T_0$ quotient, as well as the soberification, of the monad $\mathbb{U}$.
\end{rem}

\section{\bf Other examples}
We consider two examples. The first one shows that the condition $LRT\cong T$ in Theorem \ref{thm: existence of a left adjoint when RLT=T} is sufficient but not necessary for the categories of algebras \powcat{\cat{C}}{T} and \powcat{\cat{B}}{M} to be equivalent. The second example may provide a way of ``inducing'' constructions such as compactifications or completions (\cite{Hol00b}) to coreflective subcategories. 
\subsection{\bf {$\sigma$-frames}}
Consider the category \cat{\sigma Frm} of $\sigma$-frames (\cite{BanMat2003,Mad1991}), i.e. distributive lattices $L$ satisfying 
$$
\meet{a}{\left(\ds{\bigvee} S\right)}=\ds{\bigvee\set{\meet{a}{s}\ |\ s\in S}}
$$
for all $a\in L$ and any countable $S\subseteq L$, with $\sigma$-frame homomorphisms, which are lattice homomorphisms preserving countable joins. We have a composition of adjunctions
\begin{center}
\compAdj{\cat{DLat}}{\mathcal{H}_{\sigma}}{U_{\sigma}}{\cat{\sigma Frm}}{\mathcal{H}}{U}{\cat{Frm},}
\end{center}
where $U_\sigma, U$ are faithful forgetful functors and $\mathcal{H}_{\sigma},\mathcal{H}$ are left adjoint free functors (\cite{Mad1991}). We have $\comp{\comp{\mathcal{H}}{\mathcal{H}_{\sigma}}}{\comp{U_{\sigma}}{U}}=\id$. Also $\comp{\mathcal{H}_{\sigma}}{U_{\sigma}}=\id_{\sigma}$ itself is an endofunctor underlying a comonad $\mathbb{K}_{\sigma}=(\id_{\sigma},c_{\sigma},\gamma_{\sigma})$ on \cat{\sigma Frm}. Therefore
$$
\comp{\mathcal{H}}{\comp{\id_{\sigma}}{U}}\cong \id.
$$
Here $\id_{\sigma}L$ takes all countably generated ideals on $L$. We recall the following known result.
\begin{thm}
\label{thm: coalgebras of countably generated ideals}
The coalgebras of $\id_{\sigma}$ are precisely the stably continuous $\sigma$-frames or compact regular $\sigma$-biframes, together with $\sigma$-proper $\sigma$-frame homomorphisms or $\sigma$-biframe homomorphisms (\cite{BanMat2003,Mat2001}).
\end{thm}
\noindent
We denote the above category by \cat{StK\sigma Frm}. 
\begin{prop}
\label{prop: adjunction between stably continuous frames and stabl continuous sigma-frames}
The comparison functor \map{K}{\cat{StK\sigma Frm}}{\cat{StKFrm}} admits a right adjoint $K^*$.
\end{prop}
\begin{proof}
\cat{\sigma Frm} is complete. As an algebraic category\footnote{This is because \cat{\sigma Frm} is equationally presentable and the forgetful functor to \cat{DLat}, which is monadic over the category of sets and functions, admits a left adjoint. Such an argument can be seen in \cite[Section II.1.2]{Joh82}.}, \cat{\sigma Frm} admits an $(\mathcal{E},\mathcal{M})$-factorisation system, where $\mathcal{E}$ is the class of surjections (or regular epimorphisms) and $\mathcal{M}$ the class of injections (or monomorphisms) (\cite[Proposition I.3.8]{Joh82}). The functor $\id_{\sigma}$ preserves injections (\cite{Mad1991}) and \cat{\sigma Frm} is $\mathcal{M}$-well-powered. The result follows from Theorem \ref{thm: reflexive coequalizers imply existence of left adjoint} and Theorem \ref{thm: cocompleteness of algebras with (E,M) fact system}.
\end{proof}
\noindent
This adjunction is unlikely to be an equivalence. However, we know from \cite[Proposition 2.4]{BanMat2003} that the two categories are equivalent, albeit with a functor that is qualitatively different from the restriction $K^*$ of $U$. 
\begin{prop}
\label{prop: LRT not equal to T but there is equivalence}
The categories in Proposition \ref{prop: adjunction between stably continuous frames and stabl continuous sigma-frames} are equivalent, but $\comp{(\comp{U}{\mathcal{H}})}{\id_{\sigma}}\neq\id_{\sigma}$.
\end{prop}
\begin{proof}
Equivalence is shown in \cite[Proposition 2.4]{BanMat2003}. Now, since we have $\comp{(\comp{U}{\mathcal{H}})}{\id_{\sigma}}\cdot U=\comp{U}{\id}$, it is enough to show that there is a frame $L$ such that $U(\id(L))\neq\id_{\sigma}(U(L))$. But this equates to saying that there are ideals on $L$ which are not countably generated.
\end{proof}
\noindent
The proposition shows that the condition $LRT\cong T$ in Theorem \ref{thm: existence of a left adjoint when RLT=T} is not necessary.

\subsection{\bf Compactification of discrete spaces}
It is known that the ultrafilter monad $(\beta,\mu,\eta)$ on the category \cat{Set} of sets and functions (\cite{Man} and \cite[Theorem III.2.4]{Joh82}) can be ``lifted'' topologically to a monad $\mathbb{U}=(\mathfrak{U},\mu',\eta')$ on the category \cat{Top} (\cite{Low,Sal00,Raz}). In the adjunction
$$
\Adj{\cat{Set}}{F}{U}{\cat{Top}}
$$
where $U$ is the forgetful functor, and $F$ the {\em discrete functor}, we have $\comp{U}{\comp{\mathfrak{U}}{F}}=\beta$. As $F$ is a full embedding, \cat{Set} is equivalent to the subcategory of discrete spaces and continuous functions, which is then coreflective in \cat{Top}. Thus, compactifications (and completions) on \cat{Top} can be induced on its coreflective subcategories. 

\noindent
The algebras \powcat{\cat{Top}}{U} are described in \cite[Lemma 4.12 and Theorem 4.15]{BezGabJib2015} and \cite[Section III.5.6]{Low}. They are equivalent to the pre-ordered compact Hausdorff spaces and ordered continuous maps. It is known that \cat{KHaus} is the category of algebras of $(\beta,\mu,\eta)$ (\cite{Man}), and that it is reflective inside \powcat{\cat{Top}}{U} (\cite{Low}). However, it is also coreflective in \powcat{\cat{Top}}{U}, as we will show here. This is similar to the relationship between \cat{KHaus} and \cat{StKSp} (\cite[Lemma 7 and Proposition 37]{Raz}).
\begin{prop}
\label{prop: compact Hausdorff spaces are coreflective in top}
\cat{KHaus} is coreflective inside \powcat{\cat{Top}}{U}.
\end{prop}
\begin{proof}
\cat{Top} is cocomplete and admits an $(\mathcal{E},\mathcal{M})$-factorisation system, where $\mathcal{E}$ is the class of all surjective continuous maps and $\mathcal{M}$ the class of embeddings, which are extremal monomorphisms. The functor $\mathfrak{U}$ preserves $\mathcal{E}$ (\cite[Section 5.1]{Raz}). \cat{Top} is $\mathcal{E}$-well-copowered as a concrete category over \cat{Set}. Therefore the comparison functor \map{K}{\powcat{\cat{Top}}{U}}{\cat{KHaus}} admits a left adjoint which a restriction of the embedding $F$.
\end{proof}

\begin{rem}
\label{rem: Vietoris monad}
In \cite[Section 2]{HofNor2014}, the Vietoris monad $\mathbb{V}=(V,m,e)$ on \cat{StKSp} is ``transferred'' to \cat{KHaus} through the forgetful and coreflection functors between the two subcategories of \cat{Top}. This gives a monad $\mathbb{\hat{V}}=(\hat{V},m,e)$ on \cat{KHaus} where $\hat{V}X$ is the patch space of $VX$. This particular example illustrates the benefits of inducing a monad to a coreflective subcategory; this seems to be of importance for duality theories.   
\end{rem}

%%------------------------------------------------------------------------------------------------------------------
\bibliographystyle{plain}

%----------------------------------------------------------------------------------------

\end{document}